\font\script=eusm10.
\font\sets=msbm10.
\font\stampatello=cmcsc10.

\def\1{{\bf 1}}
\def\square{\hbox{\vrule\vbox{\hrule\phantom{s}\hrule}\vrule}}
\def\defineq{\buildrel{def}\over{=}}
\def\defin{\buildrel{def}\over{\Longleftrightarrow}}
\def\C{\hbox{\sets C}}
\def\N{\hbox{\sets N}}
\def\Corr{\hbox{\script C}}
\def\divisor{\hbox{\bf d}}
\def\sporadic{\sum_{{x<n\le x+h}\atop {n\equiv a\bmod q}}f(n)}
\def\short{\sum_{x<n\le x+h}f(n)}

\par
\centerline{{\bf Averages of short correlations: a note}\footnote{}{MSC 2010 : $11N37$, $11B25$, $11N36$ --- {\it Keywords}: correlations, short intervals, arithmetic progressions}}
\smallskip
\par
\centerline{\stampatello giovanni coppola}

\bigskip
\bigskip
\bigskip

\par
\centerline{\stampatello 1. Introduction and statement of main results.}
\smallskip
\par
\noindent
We have studied the so-called \lq \lq {\stampatello correlations}\rq \rq (mainly, in particular case $g=f$ of {\it autocorrelations}) of two arithmetic functions $f,g:\N \rightarrow \C$, that are also called (much more, in the modular forms environment) \lq \lq {\it shifted convolution sums}\rq \rq, i.e. (in this paper we do not conjugate $g(n-a)$, even if $g$ values are complex) 
$$
\Corr_{f,g}(N,a)\defineq \sum_{n\sim N}f(n)g(n-a), 
$$
\par
\noindent
where $a$ is called the {\stampatello shift}, the abbreviation $n\sim N$ in sums stands for $N<n\le 2N$, for integers $N>0$ (more generally, $n\sim X$ is $X<n\le 2X$, for real $X>0$). Actually, our \lq \lq Generations\rq \rq, to abbreviate [CL1] title, displays, with different approaches for the asymptotic and upper bound estimates, many different types of averages, for such quantities; mainly, it deals with three generations, namely three kinds of averages starting with the easiest, i.e., what we call the {\stampatello first generation} of correlation averages: 
$$
\sum_{a\le H}\Corr_f(a), 
$$
\par
\noindent
where we abbreviate $a\le H$ for $1\le a\le H$ and $\Corr_f(a)$ for $\Corr_{f,f}(N,a)$, in above notation (in fact, an autocorrelation). Hereafter we stick to {\it short interval} averages, namely $H=o(N)$, when $N\to \infty$. In the following, $h\ll HN^{-\varepsilon}$ will appear, as a kind of length for \lq \lq super-short intervals\rq \rq. (Compare our main results.) Hereafter the Vinogradov notation $\ll$ and the synonymous Landau's $O$ notation will be used; also, with subscripts to indicate the $O-$constant dependence. (Typically, $\ll_{\varepsilon}$ depends on arbitrarily small $\varepsilon>0$.) 

\medskip

\par
In our [CL1] study and in the subsequent [CL2] it is clear that the correlations of shift $a$ are immediately linked (see [E], for a great exposition) to the arithmetic progressions, with residue class $a$ : 
$$
\Corr_{f,g}(N,a)=\sum_{q\le Q}g'(q)\sum_{{n\sim N}\atop {n\equiv a\bmod q}}f(n), 
$$
\par
\noindent
assuming that the so-called {\stampatello Eratosthenes transform} of our $g$, $g'\defineq g\ast \mu$ (with $\ast$ the Dirichlet product and $\mu$ M\"{o}bius arithmetic function), has support up to $Q$ (see the following). Since we clearly assumed to have $Q\ll N$, there's no problem with the number of such $n$, in the inner sum (in worst case $\gg 1$ of them).
\par
In this paper we study the (a priori) much harder problem of estimating sums of the kind (here $H\ge h$ and $H=o(x)$, as we may assume $N\ll x\ll N$): 
$$
\sum_{a\le H}\sum_{x<n\le x+h}f(n)g(n-a)=\sum_{a\le H}\sum_{q\le Q}g'(q)\sporadic 
$$
\par
\noindent
where the \lq \lq short\rq \rq \enspace length $h\le H$, instead of previous \lq \lq long\rq \rq \enspace length $N$, renders this kind of {\it first  generation for short correlations} much more difficult (at least \lq \lq for arithmetic progressions\rq \rq, compare Theorem 0), as 
$$
q>h 
\enspace \Longrightarrow \enspace
\sporadic = \cases{f(n_a) & if $n_a\in (x,x+h]$\cr 0 & otherwise}
$$
\par
\noindent
(see the proof of Lemma in $\S2$) is involving \thinspace \lq \lq {\it sporadic sums}\rq \rq, \thinspace i.e. sums with \lq \lq {\it at most one term and}\rq \rq, in a probabilistic sense, \lq \lq {\it vanishing most of the time}\rq \rq. (In the above example, see Lemma proof, this simply means : the number of $a$s for which the $n-$sum doesn't vanish turns out to be \lq \lq small\rq \rq, compared to all $a$s.) 
\par				
However, an elementary argument (our Lemma) shows that we may gather $H\ge h$ residue classes, so to avoid 
sporadically appearing terms !
Here, the fact that we consider a kind of average, instead of single correlations, is vital ! (Actually, without any kind of average, short correlations may not be treated at all !) 

\bigskip

\par
Before going on, we need the following two definitions (for which, compare [CL1] \& [CL2]). 

\medskip

\par
\noindent
We say, for a general arithmetic function $f:\N \rightarrow \C$, 
$$
f \enspace \hbox{\rm is} \enspace \hbox{\stampatello essentially \thinspace bounded} 
\enspace \defin \enspace 
\forall \varepsilon>0, \enspace f(n)\ll_{\varepsilon} n^{\varepsilon}, 
\enspace \hbox{\rm as} \enspace n\to \infty. 
$$
\par
\noindent
(Notice : compare this with Ramanujan's Conjecture in Selberg Class and in modular forms environment.) 

\smallskip

\par
We say $f:\N \rightarrow \C$ is a {\stampatello sieve function of range} $Q$ when the Eratosthenes transform of $f$, say, $f':\N \rightarrow \C$ (so that $f=f'\ast \1$ by M\"{o}bius inversion : hereafter $\1$ is the constant-$1$ arithmetic function), is supported up to $Q$ and essentially bounded; in other words, the sieve functions of range $Q$ are representable as 
$$
f(n)\defineq \sum_{q|n,q\le Q}f'(q), 
$$
\par
\noindent
for a certain arithmetic function $f'$ (the $f$ Eratosthenes transform), satisfying \enspace $f'(q)\ll_{\varepsilon} q^{\varepsilon}$. Since, by M\"{o}bius inversion, \lq \lq $f$ is essentially bounded\rq \rq \enspace is equivalent to \lq \lq $f\ast \mu$ is essentially bounded\rq \rq, \enspace we may also identify sieve functions of range $Q$ with truncated divisor sums up to $Q$ that satisfy Ramanujan Conjecture. 

\medskip

\par
We give, for first, our easiest result for averages of short correlations. 
\smallskip
\par
\noindent {\bf Theorem 0}. {\it Let } $h,H,x\in \N$ {\it with } $h\le H$ {\it and } $H=o(x)$, {\it as } $x\to \infty$. {\it Take any } {\stampatello essentially bounded} {\it arithmetic functions } $f,g:\N \rightarrow \C$. {\it Then}
$$
\sum_{a\le H}\sum_{x<n\le x+h}f(n)g(n-a)
=\left(\sum_{a\le H}g(x-a)\right)\left(\sum_{x<n\le x+h}f(n)\right)+O_{\varepsilon}\left(x^{\varepsilon} h^2\right). 
$$
\smallskip
\par
\noindent
{\bf Proof} is exhibited here, as it is $1-$line (setting $b:=a-n+x$ and back $a$ instead of $b$, using $n-x\ll h$) : 
$$
\sum_{a\le H}\sum_{x<n\le x+h}f(n)g(n-a)
=\sum_{x<n\le x+h}f(n)\sum_{a\le H}g(x-(a-n+x))
=\sum_{x<n\le x+h}f(n)\sum_{a\le H}g(x-a)+O_{\varepsilon}\left(x^{\varepsilon} h^2\right). 
$$
\hfill $\square$ 
\par
\centerline{(A symbol \enspace \thinspace $\square$ \thinspace \enspace will signify, as usual, the end of a Proof, likewise \enspace \thinspace  $\diamond$ \thinspace \enspace the end of a Remark.)}
\smallskip
\par
\noindent {\bf Remark 0}. The remainders are non-trivial when \enspace $h\ll H/x^{\varepsilon}$, i.e., the $H-$average is \lq \lq long enough\rq \rq. 
\hfill $\diamond$ 
\smallskip
\par
\noindent {\bf Remark 1}. See that the following results for short correlations averages are intended both to give some further insight and to get new results, from comparing with Theorem 0; which is, of course, far simpler ! 
\hfill $\diamond$ 

\medskip

\par
Our second main result for averages of short correlations will be proved in $\S3$. 
\smallskip
\par
\noindent {\bf Theorem 1}. {\it Let } $h,H,Q,x\in \N$ {\it with } $h\le H<Q$ {\it and } $Q\ll x$, $H=o(x)$ {\it as } $x\to \infty$. {\it Take any } {\stampatello essentially bounded} $f:\N \rightarrow \C$ {\it and a} {\stampatello sieve function } $g:\N \rightarrow \C$ {\stampatello of range} $Q$ {\it with Eratosthenes transform } $g'$. {\it Then}
$$
\sum_{a\le H}\sum_{x<n\le x+h}f(n)g(n-a)
=\left(H\sum_{q\le h}{{g'(q)}\over q} + \sum_{h<q\le H}g'(q)\left[{H\over q}\right] 
                                       + \sum_{{h<q\le Q}\atop {\left\{{x\over q}\right\}\le \left\{{H\over q}\right\}}}g'(q)
  \right)\short + 
$$
$$
 +O_{\varepsilon}\left(x^{\varepsilon} h^2\right). 
$$ 
\par
\noindent
See, for remarks and comments, soon after its proof. 

\medskip

\par				
Rather an identical argument also gives the following companion result, proved in $\S3$, too. 
\smallskip
\par
\noindent {\bf Theorem 2}. {\it Let the same hypotheses of Theorem 1 hold.} {\it Then}
$$
\sum_{a\le H}\sum_{x<n\le x+h}f(n)g(n-a)=\left(\sum_{a\le H}g(x-a)\right)\left(\short\right) + O_{\varepsilon}\left(x^{\varepsilon}h^2\right) + 
$$
$$
+ \left(H\sum_{q\le h}{{g'(q)}\over q}\right)\short + \sum_{a\le H}\sum_{h<q\le H}g'(q)\sporadic 
 - \left(\sum_{a\le H}\sum_{{q\le H}\atop {q|x-a}}g'(q)\right)\short. 
$$

\medskip

\par
In passing, from the comparison of Theorem 0 and Theorem 2, the following result immediately follows. 
\smallskip
\par
\noindent {\bf Corollary}. {\it Let } $h,H,x\in \N$ {\it with } $h\le H$ {\it and } $H=o(x)$, {\it as } $x\to \infty$. {\it Take any } {\stampatello essentially bounded} {\it arithmetic functions } $f,g':\N \rightarrow \C$. {\it Then}

$$
\left(H\sum_{q\le h}{{g'(q)}\over q}\right)\short
 + \sum_{a\le H}\sum_{h<q\le H}g'(q)\sporadic 
  - \left(\sum_{a\le H}\sum_{{q\le H}\atop {q|x-a}}g'(q)\right)\short
\ll_{\varepsilon} x^{\varepsilon}h^2. 
$$
\par
\noindent
{\it This entails, choosing first a vanishing $g'$ in the interval } $[1,h]$ {\it and, then, vanishing in } $(h,H]$, {\stampatello both} 
$$
\sum_{a\le H}\sum_{h<q\le H}g'(q)\sporadic = \left(\sum_{a\le H}\sum_{{h<q\le H}\atop {q|x-a}}g'(q)\right)\short 
 +O_{\varepsilon}\left(x^{\varepsilon}h^2\right) 
$$
\par
\noindent
{\stampatello and} 
$$
\left(\sum_{a\le H}\sum_{{q\le h}\atop {q|x-a}}g'(q)\right)\short = \left(H\sum_{q\le h}{{g'(q)}\over q}\right)\short 
 +O_{\varepsilon}\left(x^{\varepsilon}h^2\right). 
$$
\smallskip
\par
\noindent {\bf Remark 2}. The further hypotheses on $g$ and $Q$, from Theorem 2 (i.e. from Th.m 1), are not used at all ! \hfill $\diamond$ 

\medskip

\par
As we saw before, it is then immediate to glue together these {\it super-short intervals} to get, as a corollary to Theorem 0, a result for the first generation of correlations (these are the \lq \lq long ones\rq \rq, this time). 
\smallskip
\par
\noindent {\bf Theorem 3}. {\it Let } $H,D,Q,N\in \N$ {\it with } $H<Q$ {\it and } $D,Q\ll N$, $H=o(N)$ {\it as } $N\to \infty$. {\it Choose an integer } $h$ {\it such that } $h\to \infty$ {\it as } $N\to \infty$. {\it Take two} {\stampatello sieve functions} $f,g:\N \rightarrow \C$ {\stampatello of ranges } $D,Q$, {\it resp. } {\it Then} 
$$
\sum_{a\le H}\sum_{n\sim N}f(n)g(n-a) = \sum_{a\le H}\sum_{n\sim N}f(n)g(h[n/h]-a) + O_{\varepsilon}\left(N^{\varepsilon}H\left({{Nh}\over H}+{N\over h}+h\right)\right). 
$$
\par
\noindent
{\it Furthermore, averaging up to such a fixed integer,} 
$$
\sum_{a\le H}\sum_{n\sim N}f(n)g(n-a) = \sum_{a\le H}\sum_{n\sim N}f(n){1\over h}\sum_{m\le h}g(m[n/m]-a) + O_{\varepsilon}\left(N^{\varepsilon}H\left({{Nh}\over H}+{N\over h}+h\right)\right). 
$$
\par
\noindent {\bf Remark 3}. The second formula follows by estimating trivially for $m\le L:=\log N \to \infty$ and for \thinspace $L<m\le h$ \thinspace applying first formula (all \lq \lq $\log$s\rq \rq \thinspace and so on are inside $N^{\varepsilon}$), whose proof we leave to the interested reader. \hfill $\diamond$ 
\medskip
\par
We'll study applications of Theorem 3 in future papers (and versions of present one). 
\medskip
\par
\noindent
The paper is organized as follows: our Lemma is given and proved in $\S2$, then $\S3$ proves main results, namely Theorem 1 and Theorem 2. 
(Notice that we call these \lq \lq main\rq \rq, as opposed to side results, i.e., the Lemma and the Remarks. 
These, in turn, have their own interest, not only in view of proving Theorems.) 

\bigskip


\par				
\centerline{\stampatello 2. Pinch Lemma.}
\smallskip
\par
\noindent
We start with so-called \lq \lq {\stampatello Pinch Lemma}\rq \rq, since it pinches, so to speak, integers in (very) short intervals, \lq \lq squeezing\rq \rq \enspace them to the left interval extreme (see Remark 4, after the proof). In fact, these short intervals contain at most one integer in the specified residue class modulo $q$, since $q>h$, where $h$ is the interval's length. We actually rebuild the whole interval, after summing over $H\ge h$ residue classes. 
\smallskip
\par
\noindent {\bf Lemma}. {\it Let } $h,H,Q\in \N$ {\it with } $h\le H<Q$ {\it and } $Q\ll x$, $H=o(x)$, {\it as } $x\to \infty$. {\it Take two} {\stampatello essentially bounded} {\it arithmetic functions } $f,g':\N \rightarrow \C$. {\it Then} {\stampatello both} 
$$
\sum_{a\le H}\sum_{h<q\le H}g'(q)\sporadic = \left(\sum_{h<q\le H}g'(q)\left[{H\over q}\right]+\sum_{{h<q\le H}\atop {\left\{ {x\over q}\right\}\le \left\{ {H\over q}\right\}}}g'(q)\right)\short + O_{\varepsilon}\left(x^{\varepsilon}h^2\right)
$$
\par
\noindent
{\stampatello and}
$$
\sum_{a\le H}\sum_{H<q\le Q}g'(q)\sporadic = \left(\sum_{{H<q\le Q}\atop {\left\{{x\over q}\right\}\le {H\over q}}}g'(q)\right)\short + O_{\varepsilon}\left(x^{\varepsilon}h^2\right). 
$$

\medskip

\par
\noindent
Notice we have non-triviality for the remainders when $h$ is significantly smaller than $H$, say $\log h<\log H$. 

\medskip

\par
\noindent {\bf Proof}. For the first formula write LHS (left hand side) as 
$$
\sum_{h<q\le H}g'(q)\Big(\sum_{a\le q\left[{H\over q}\right]}+\sum_{q\left[{H\over q}\right]<a\le H}\Big)\sporadic = 
$$
$$
= \Big(\sum_{h<q\le H}g'(q)\left[{H\over q}\right]\Big)\short 
 + \sum_{h<q\le H}g'(q)\sum_{a\le q\left\{{H\over q}\right\}}\sporadic 
$$
\par
\noindent
(recall \thinspace $a\le A$ \thinspace is \thinspace $1\le a\le A$) and see that \enspace $0<a\le q\{H/q\}<q$ \enspace allows to write 
$$
\sum_{a\le q\left\{{H\over q}\right\}}\sporadic = \sum_{a\le q\left\{{H\over q}\right\}}\sum_{{{x-a}\over q}<m\le {{x-a+h}\over q}}f(qm+a), 
$$
\par
\noindent
whence $m-$sum is {\stampatello sporadic} (i.e., contains at most one term, as the interval has length $h/q<1$) and writing the $m-$interval as 
$$
\left({{x-a}\over q},{{x-a+h}\over q}\right]=\left(\left[{x\over q}\right]+\left\{{x\over q}\right\}-{a\over q},\left[{x\over q}\right]+\left\{{x\over q}\right\}-{a\over q}+{h\over q}\right]
$$
\par
\noindent
we see that \thinspace $m=[x/q]$ \thinspace or \thinspace $m=[x/q]+1$ (in particular, $H=o(x)\Rightarrow qm>x-a>0 \Rightarrow m>0$), therefore (observe: next $a-$intervals are disjoint) 
$$
\sum_{{{x-a}\over q}<m\le {{x-a+h}\over q}}f(qm+a)=f\left(q\left[{x\over q}\right]+a\right)\1_{q\left\{{x\over q}\right\}<a\le q\left\{{x\over q}\right\} +h} 
 + f\left(q\left[{x\over q}\right]+q+a\right)\1_{a\le q\left\{{x\over q}\right\}+h-q}
\leqno{(\ast)}
$$
\par
\noindent
gives 
$$
\sum_{a\le q\left\{{H\over q}\right\}}\sum_{{{x-a}\over q}<m\le {{x-a+h}\over q}}f(qm+a)
=\sum_{{a\le q\left\{{H\over q}\right\}}\atop {q\left\{{x\over q}\right\}<a\le q\left\{{x\over q}\right\}+h}}f\left(q\left[{x\over q}\right]+a\right) 
 + \sum_{{a\le q\left\{{H\over q}\right\}}\atop {a\le q\left\{{x\over q}\right\}+h-q}}f\left(q\left[{x\over q}\right]+q+a\right) = 
$$
$$				
=\1_{\left\{{x\over q}\right\}\le \left\{{H\over q}\right\}-{h\over q}} \short 
 + \1_{\left\{{H\over q}\right\}-{h\over q}<\left\{{x\over q}\right\}<\left\{{H\over q}\right\}} O_{\varepsilon}\left(x^{\varepsilon}h\right) 
 + \1_{\left\{{x\over q}\right\}>1-{h\over q}}O_{\varepsilon}\left(x^{\varepsilon}h\right). 
$$
\par
\noindent
Here we abbreviate \enspace $\1_{\wp}\defineq 1$, when $\wp$ is true, $\defineq 0$ otherwise (while $\1$ is the constant-$1$ arithmetic function). 
\par
\noindent
Once summed over $h<q\le H$ with $g'(q)$, it implies first formula, because: 
$$
\sum_{h<q\le H}g'(q)\sum_{a\le q\left\{{H\over q}\right\}}\sporadic = \Big(\sum_{{h<q\le H}\atop {\left\{{x\over q}\right\}\le \left\{{H\over q}\right\}-{h\over q}}}g'(q)\Big)\short + 
$$
$$
+ O_{\varepsilon}\Big(x^{\varepsilon}\sum_{{h<q\le H}\atop {\left\{{H\over q}\right\}-{h\over q}<\left\{{x\over q}\right\}<\left\{{H\over q}\right\}}}h\Big) + O_{\varepsilon}\Big(x^{\varepsilon}\sum_{{h<q\le H}\atop {1-{h\over q}<\left\{{x\over q}\right\}<1}}h
\Big) 
=\Big(\sum_{{h<q\le H}\atop {\left\{{x\over q}\right\}\le \left\{{H\over q}\right\}}}g'(q)\Big)\short 
+ O_{\varepsilon}\left(x^{\varepsilon}h^2\right), 
$$
\par
\noindent
since the {\stampatello divisor function} is essentially bounded : $\divisor(n)\defineq \sum_{d|n}1\ll_{\varepsilon} n^{\varepsilon}$, a well-known fact we use inside {\stampatello both} 
$$
\sum_{{h<q\le H}\atop {\left\{{H\over q}\right\}-{h\over q}<\left\{{x\over q}\right\}\le \left\{{H\over q}\right\}}}1
=\sum_{{h<q\le H}\atop {q\left\{{H\over q}\right\}-h<q\left\{{x\over q}\right\}\le q\left\{{H\over q}\right\}}}1
=\sum_{h<q\le H}\sum_{{0\le a<h}\atop {a=q\left\{{H\over q}\right\}-q\left\{{x\over q}\right\}}}1
=\sum_{h<q\le H}\sum_{{0\le a<h}\atop {x+a\equiv H\bmod q}}1= 
$$
$$
=\sum_{0\le a<h}\sum_{{h<q\le H}\atop {q|x+a-H}}1
\le \sum_{0\le a<h}\divisor(x+a-H)
\ll_{\varepsilon} x^{\varepsilon} h 
$$
\par
\noindent
(here and in the following we use that $a-$intervals have length $\le h<q$) {\stampatello and} 
$$
\sum_{{h<q\le H}\atop {1-{h\over q}<\left\{{x\over q}\right\}<1}}1
=\sum_{{h<q\le H}\atop {q-h<q\left\{{x\over q}\right\}<q}}1
=\sum_{h<q\le H}\sum_{{0<a<h}\atop {a=q-q\left\{{x\over q}\right\}}}1
=\sum_{0<a<h}\sum_{{h<q\le H}\atop {q|x+a}}1
\le \sum_{0<a<h}\divisor(x+a)
\ll_{\varepsilon} x^{\varepsilon} h. 
$$
\par
\noindent
(For these and following bounds, it is {\sl vital} that the \enspace $\divisor(n)$ \enspace has \thinspace $n>0$ : thanks to $H=o(x)$ hypothesis.) 
\par
\noindent
Analogously, for the second formula, we apply $(\ast)$ to get (see that now $a<q$ is for free from $H<q$) 
$$
\sum_{a\le H}\sum_{H<q\le Q}g'(q)\sporadic = \sum_{H<q\le Q}g'(q)\sum_{a\le H}\sum_{{{x-a}\over q}<m\le {{x-a+h}\over q}}f(qm+a)= 
$$
$$
= \sum_{H<q\le Q}g'(q)\sum_{{a\le H}\atop {q\left\{{x\over q}\right\}<a\le q\left\{{x\over q}\right\}+h}}f(qm+a) 
 + O_{\varepsilon}\Big(x^{\varepsilon}\sum_{{H<q\le Q}\atop {1-{h\over q}<\left\{{x\over q}\right\}<1}}h\Big)= 
$$
$$
= \sum_{{H<q\le Q}\atop {q\left\{{x\over q}\right\}\le H-h}}g'(q)\short 
 + O_{\varepsilon}\Big(x^{\varepsilon}h\sum_{{H<q\le Q}\atop {H-h<q\left\{{x\over q}\right\}<H}}1\Big)
  + O_{\varepsilon}\left(x^{\varepsilon}h^2\right)= 
$$
$$
= \sum_{{H<q\le Q}\atop {q\left\{{x\over q}\right\}\le H}}g'(q)\short 
 + O_{\varepsilon}\Big(x^{\varepsilon}h\sum_{{H<q\le Q}\atop {H-h<q\left\{{x\over q}\right\}\le H}}1\Big)
  + O_{\varepsilon}\left(x^{\varepsilon}h^2\right), 
$$
\par
\noindent
similarly as above for \enspace $q-h<q\{x/q\}<q$ \enspace and we conclude by the analogous (set $a=H-q\{x/q\}$ now) 
$$
\sum_{{H<q\le Q}\atop {H-h<q\left\{{x\over q}\right\}\le H}}1=\sum_{0\le a<h}\sum_{{H<q\le Q}\atop {q|x+a-H}}1
\le \sum_{0\le a<h}\divisor(x+a-H)
\ll_{\varepsilon} x^{\varepsilon} h.
$$
\hfill $\square$
\par
\noindent {\bf Remark 4}. We may write Lemma's second formula as (we pinch $n\equiv a\bmod q$ so it \lq \lq becomes\rq \rq \thinspace $x\equiv a\bmod q$) 
$$
\sum_{a\le H}\sum_{H<q\le Q}g'(q)\sporadic = \left(\sum_{a\le H}\sum_{{H<q\le Q}\atop {q|x-a}}g'(q)\right)\short + O_{\varepsilon}\left(x^{\varepsilon}h^2\right). 
$$
\hfill $\diamond$ 

\vfill
\eject 

\par				
\centerline{\stampatello 3. Proof of main results.}
\smallskip
\par
\noindent{\bf Proof (Th.1)}. Taking $g'$ in the Pinch Lemma to be the Eratosthenes Transform of $g$, open \enspace $g(n-a)$ : 
$$
\sum_{a\le H}\sum_{x<n\le x+h}f(n)g(n-a)=\sum_{a\le H}\sum_{q\le Q}g'(q)\sporadic 
= \sum_{a\le H}\left(\sum_{q\le h}+\sum_{h<q\le H}+\sum_{H<q\le Q}\right)g'(q)\sporadic 
$$
$$
=I+II+III, 
$$
\par
\noindent
say, where 
$$
I:=\sum_{a\le H}\sum_{q\le h}g'(q)\sporadic 
= \sum_{q\le h}g'(q)\left(\sum_{a\le q\left[{H\over q}\right]}+\sum_{q\left[{H\over q}\right]<a\le H}\right)\sporadic = 
$$
$$
= \sum_{q\le h}g'(q)\left[{H\over q}\right]\short + O_{\varepsilon}\left(x^{\varepsilon} h\sum_{x<n\le x+h}1\right) 
= H\sum_{q\le h}{{g'(q)}\over q}\short + O_{\varepsilon}\left(x^{\varepsilon} h\sum_{x<n\le x+h}1\right)= 
$$
$$
= \left(H\sum_{q\le h}{{g'(q)}\over q}\right)\short + O_{\varepsilon}\left(x^{\varepsilon} h^2\right), 
$$
\par
\noindent
while by the Pinch Lemma 
$$
II:=\sum_{a\le H}\sum_{h<q\le H}g'(q)\sporadic 
=\Big(\sum_{h<q\le H}g'(q)\left[{H\over q}\right] + \sum_{{h<q\le H}\atop {\left\{{x\over q}\right\}\le \left\{{H\over q}\right\}}}g'(q)
  \Big)\short 
 +O_{\varepsilon}\left(x^{\varepsilon} h^2\right) 
$$
\par
\noindent
and 
$$
III:=\sum_{a\le H}\sum_{H<q\le Q}g'(q)\sporadic 
=\Big(\sum_{{H<q\le Q}\atop {\left\{{x\over q}\right\}\le {H\over q}}}g'(q)
  \Big)\short 
 +O_{\varepsilon}\left(x^{\varepsilon} h^2\right)= 
$$
$$
=\Big(\sum_{{H<q\le Q}\atop {\left\{{x\over q}\right\}\le \left\{{H\over q}\right\}}}g'(q)
  \Big)\short 
 +O_{\varepsilon}\left(x^{\varepsilon} h^2\right). 
$$
\hfill $\square$ 

\medskip

\par
\noindent {\bf Remark 5}. However, see the definition of $II$ in the Theorem proof, same hypotheses give 
$$
\sum_{a\le H}\sum_{x<n\le x+h}f(n)g(n-a)
= \left(H\sum_{q\le h}{{g'(q)}\over q} + \sum_{a\le H}\sum_{{H<q\le Q}\atop {q|x-a}}g'(q)\right)\short + 
$$
$$
+ \sum_{a\le H}\sum_{h<q\le H}g'(q)\sporadic + O_{\varepsilon}\left(x^{\varepsilon}h^2\right), 
$$
\par
\noindent
thanks to Remark 4 (after Pinch Lemma proof), too. 
\hfill $\diamond$ 

\medskip

\par
\noindent{\bf Proof (Th.2)}. Follows from previous Remark and 
$$
\sum_{{H<q\le Q}\atop {q|x-a}}g'(q)
=\sum_{{q\le Q}\atop {q|x-a}}g'(q)-\sum_{{q\le H}\atop {q|x-a}}g'(q)
=g(x-a)-\sum_{{q\le H}\atop {q|x-a}}g'(q). 
$$
\hfill $\square$ 

\bigskip

\par				
\centerline{\stampatello References.}
\smallskip
\item{\bf [CL1]} Coppola, G. and Laporta, M. \thinspace - \thinspace {\sl Generations of correlation averages} \thinspace - \thinspace J. Numbers Volume 2014 (2014), Article ID 140840, 13 pages http://dx.doi.org/10.1155/2014/140840 (see draft arxiv:1205.1706v3 too) 
\item{\bf [CL2]} Coppola, G. and Laporta, M. - {\sl Sieve functions in arithmetic bands} - http://arxiv.org/abs/1503.07502v3 
\item{\bf [E]} Elliott, P.D.T.A. \thinspace - \thinspace {\sl On the correlation of multiplicative and the sum of additive arithmetic functions} Mem. Amer. Math. Soc. {\bf 112} (1994), no. 538, viii+88 pp. $\underline{\tt MR\thinspace 95d\!:\!11099}$ 

\vfill

\par
\noindent
{\stampatello Giovanni Coppola}
\par
\noindent
Postal address: Via Partenio 12, 
\par
\noindent
83100 Avellino (AV), ITALY
\par
\noindent
e-mail : giocop@interfree.it
\par
\noindent
wpage : www.giovannicoppola.name

\bye